\documentclass[12pt]{amsart}
\usepackage{amssymb}
\usepackage{amsmath}
\usepackage{latexsym}
\usepackage{amsthm}
\makeatletter
\def\thmhead@plain#1#2#3{%
 \thmname{#1}\thmnumber{\@ifnotempty{#1}{ }#2}%
 \thmnote{ \the\thm@notefont(#3)}}
\let\thmhead\thmhead@plain
\def\swappedhead#1#2#3{%
 \thmnumber{#2}\thmname{\@ifnotempty{#2}{. }#1}%
 \thmnote{ \the\thm@notefont(#3)}}
\makeatother

\theoremstyle{definition} 
\newtheorem{definition}{Definition}[section]
\newtheorem{remark}[definition]{Remark}

\theoremstyle{plain}      
\newtheorem{proposition}[definition]{Proposition}
\newtheorem{theorem}[definition]{Theorem}
\newtheorem{corollary}[definition]{Corollary}

\newcommand{\Var}{\operatorname{Var}}
\newcommand{\Cov}{\operatorname{Cov}}

\begin{document}
\title[Block-hierarchical covariance decompositions]{Block-hierarchical covariance decompositions for finite-block additive functionals}

\author{Abbas Alhakim}
\address{Department of Mathematics, American University of Beirut, Beirut, Lebanon}
\email{aa145@aub.edu.lb}
\thanks{ORCID: 0000-0002-2995-3740.}
\subjclass[2020]{Primary 60J10, 60F05; Secondary 47A10, 60G10, 62M15.}

\keywords{finite-block additive functionals, overlapping windows, Markov chains, Green--Kubo covariance, spectral decomposition, Hilbert-space projections. }

\date{July, 2026.}


\begin{abstract}
We study additive functionals of stationary Markov chains whose observables depend on a fixed finite block of consecutive states. Such block observables arise naturally in sliding-window statistics, pattern counts, and local dependence analysis. In the independent setting, additive functionals of overlapping finite blocks are known to have a covariance operator with integer spectrum \(0,1,\ldots,k\), and the eigenvalue-one component represents the information first detectable at block length \(k\). We study the corresponding problem when the underlying sequence is a stationary Markov chain. For a fixed block observable \(f(X_t,\ldots,X_{t+k-1})\), we introduce a Hilbert-space decomposition that separates information already contained in shorter consecutive blocks from the genuinely new block-\(k\) component, called the incremental component. We show that this component persists as  an eigenvalue-one component under Markovian dependence: on this space the Green--Kubo covariance converges trivially and the covariance operator acts as the identity. More generally, the integers \(1,\ldots,k-1\) are shown to arise hierarchically as eigenvalues of the Green--Kubo operator. In the reversible case, the remaining covariance structure is described through a boundary-corrected one-coordinate marginal spectrum determined by the base transition operator. Reversible two-state chains and Gaussian AR(1) models illustrate the theory through concrete spectral formulas.
\end{abstract}

\maketitle 

\section{Introduction}

Let \((X_t)_{t\in\mathbb Z}\) be a stationary Markov chain.  A broad class of
statistics and observables of such a process is formed by evaluating a function on
a finite block of consecutive states and summing along the path
\[
        \sum_{t=1}^n f(X_t,\ldots,X_{t+k-1}).
\]
We call such an \(f\) a finite-block observable.  Equivalently, it is a cylinder
observable on the stationary path space, depending only on \(k\) consecutive
coordinates.  For fixed \(k\), the lifted process
\[
        Y_t^{(k)}=(X_t,\ldots,X_{t+k-1})
\]
is itself a Markov chain, and the preceding sum is an ordinary additive functional
of this \(k\)-window chain.

Classical limit theory for additive functionals of Markov chains, through the
Green--Kubo formula, Poisson equations, and martingale approximation, provides
powerful tools for studying the asymptotic behavior of such sums; see, for example,
Gordin \cite{Gordin1969}, Kipnis and Varadhan \cite{KipnisVaradhan1986}, and Wu
\cite{Wu2005}. Applied to the lifted chain \(Y_t^{(k)}\), this theory treats \(f(X_t,\ldots,X_{t+k-1})\) as an ordinary state function on \(E^k\). It therefore does not address a different, block-hierarchical question: what part of \(f(X_t,\ldots,X_{t+k-1})\) is genuinely new at block length \(k\), and what part is already representable using shorter consecutive blocks?

This question is especially visible in sliding-window statistics. Increasing the
window length is often intended to capture additional local dependence, ordering,
or pattern information. In finite-alphabet and symbolic settings, the related
question of whether longer blocks carry genuinely new information is traditionally
studied through likelihood-ratio methods, information criteria, or context-tree
estimators, which compare conditional probability models of different memory
lengths; see, for example,
\cite{Tong1975,Katz1981,BuhlmannWyner1999,CsiszarTalata2006}. The viewpoint
developed here is different: it does not estimate the conditional law, but instead
identifies Hilbert-space directions of \(k\)-block observables that are not
reducible to adjacent shorter blocks. Such sliding-window constructions appear in
overlapping serial tests of randomness, approximate entropy, pattern-count
statistics, and local dependence analysis; see, for example, Good~\cite{Good1953},
Pincus~\cite{Pincus1991}, Alhakim~\cite{Alhakim2004}, and Alhakim, Kawczak, and
Molchanov~\cite{AlhakimKawczakMolchanov2004}. In the independent setting, the associated
overlapping-window chain has a nilpotent covariance structure and its covariance
operator has the integer spectrum
\[
        \{0,1,\ldots,k\}.
\]
The eigenvalue-one component is the part first detectable at block length \(k\).

The present paper studies what remains of this structure under Markovian dependence.
Once the underlying sequence itself has memory, the nilpotent structure of the
independent case disappears, and the covariance operator is governed by the
transition dynamics of the base chain through the Green--Kubo formula.  Nevertheless,
a substantial part of the overlap hierarchy persists.  We show that the incremental
block-\(k\) component has an intrinsic Hilbert-space definition as the orthogonal
complement of the information contained in shorter consecutive blocks.  For a
stationary Markov chain this space admits the dynamic characterization
\[
        \mathcal I_k=\ker P_k\cap\ker P_k^*,
\]
where \(P_k\) is the transition operator of the \(k\)-window chain.  Consequently, the Green--Kubo series is trivial on \(\mathcal I_k\), and \(\mathcal I_k\)
is an eigenvalue-one subspace of the covariance operator.

Beyond this universal component, the Markovian covariance structure separates into
two further mechanisms.  First, incremental components of smaller block lengths give
rise hierarchically to the integer eigenvalues \(1,\ldots,k-1\).  Second, in the
reversible case, the one-coordinate marginal component is described through the
spectral resolution of the transition operator of the base chain.  Thus the independent integer
spectrum is partly preserved by the overlap geometry and partly deformed by the
Markov dynamics of the underlying chain.

The present framework belongs to a broader family of orthogonal decompositions in
probability and statistics.  Hoeffding-type decompositions for \(U\)-statistics
organize functions according to interaction order; see Hoeffding \cite{Hoeffding1948}
and, for dependent settings, Ho and Hsing \cite{HoHsing1997} and recent
Hilbert-space approaches such as \cite{IlIdrissi2025}.  Martingale and
Poisson-equation methods decompose additive functionals through filtration and
coboundary structures, while spectral methods for stationary processes decompose
variance contributions through frequency-domain representations; see Brillinger
\cite{Brillinger1975}.  The decomposition studied here is different: it is organized
by consecutive block length and by the covariance geometry induced by overlap.

The rest of the paper is organized as follows.  Section~2 recalls the covariance
operator and the independent nilpotent framework.  Section~3 states the main
results: the intrinsic incremental space, the \(k=2\) marginal-transition separation,
the hierarchical integer eigenvalues, the reversible spectrum, and the downward
decomposition.  Section~4 contains the proofs.  Section~5 develops reversible
two-state and Gaussian AR(1) examples, illustrating how the abstract decomposition
reduces to explicit spectral formulas.


\section{Background and notation}

We now fix the operator notation used throughout the paper. The point of this section is twofold: first, to formulate finite-block observables as additive functionals of a lifted Markov chain; and second, to recall the independent overlap structure that motivates the spectral questions studied below. We use standard notation and facts for Markov kernels and their associated operators;
see, for example, Chung \cite{KLChung1982}.

\subsection{Additive functionals and overlap covariance}

Let \(\{X_t\}_{t\in\mathbb Z}\) be a stationary Markov chain on a measurable
state space \((E,\mathcal E)\), with invariant distribution \(\nu\) and transition
kernel \(K\).  For a fixed block length  \(k\ge2\), define the lifted $k$-block, or $k$-window, chain
\begin{equation}\label{E:WindowChain}
Y_t^{(k)}=(X_t,X_{t+1},\ldots,X_{t+k-1}),
\end{equation}
which is itself a stationary Markov chain on
\(E^k\), whose dependence
structure is generated jointly by the overlap of consecutive blocks and the
transition dynamics of the underlying chain. We denote its invariant distribution by \(\pi_k\).  

Explicitly, \(\pi_k\)
is the law of a length-\(k\) stationary path of the base chain:
\[
\pi_k(dx_1,\ldots,dx_k)
=
\nu(dx_1)K(x_1,dx_2)K(x_2,dx_3)\cdots K(x_{k-1},dx_k).
\]
We denote by \(\mathcal H_k=L^2_0(\pi_k)\) the Hilbert space of square-integrable
functions \(f:E^k\to\mathbb R\) satisfying the centering condition
\[
\int f\,d\pi_k=0.
\]

The Markov transition kernel \(K\) induces a contraction on \(L^2(\nu)\), denoted
by the same symbol:
\begin{equation}\label{E:MarkovOperatorK}
(Kh)(x)=\int_E h(y)\,K(x,dy).
\end{equation}
Since \(\nu\) is stationary, \(K\) preserves constants and maps \(L^2_0(\nu)\) into
itself. 
Thus, by the Markov property of the underlying chain and the definition of the
\(k\)-block chain, the transition operator \(P_k\) of the $k$-block chain acts on
\(\mathcal H_k\) as
\begin{equation}\label{E:MarkovOperatorP}
(P_k f)(x_1,\ldots,x_k)
=
\int_E f(x_2,\ldots,x_k,z)\,K(x_k,dz).
\end{equation}

The adjoint \(P_k^*\) is the transition operator of the time-reversed window chain.
If \(K^*\) denotes the time-reversal of \(K\) with respect to \(\nu\), then
\[
P_k^* f(x_1,\ldots,x_k)
=
\int_E
f(z,x_1,\ldots,x_{k-1})\,K^*(x_1,dz).
\]
In the reversible case, \(K^*=K\).

For \(f\in \mathcal H_k\), we consider the normalized additive functional
\begin{equation}
S_n(f)=\frac1{\sqrt n}\sum_{t=1}^n f(Y_t^{(k)}).
\end{equation}
Under standard ergodicity and moment assumptions ensuring a central limit theorem
for additive functionals of Markov chains, one has (see \cite{Gordin1969,KipnisVaradhan1986})
\begin{equation}\label{E:CLT}
S_n(f)\Rightarrow N(0,\sigma_k^2(f)),
\end{equation}
where the asymptotic variance admits the Green--Kubo representation
\[
\sigma_k^2(f)
=
\Var_{\pi_k}(f)
+
2\sum_{h=1}^\infty
\Cov_{\pi_k}
\bigl(
f(Y_0^{(k)}),
f(Y_h^{(k)})
\bigr).
\]
Equivalently,
\[
\sigma_k^2(f)
=
\langle f,f\rangle
+
2\sum_{h=1}^\infty
\langle f,P_k^h f\rangle,
\]
whenever the series converges. We shall write this covariance form as 
\begin{equation}\label{E:GreenKubo}
 B_k = I+\sum_{h\geq 1}P_k^h+\sum_{h\geq 1}(P_k^*)^h, 
\end{equation}
so that 
\[ \sigma_k^2(f)=\langle B_k f,f\rangle . 
\] 
Throughout the paper, whenever \(B_k\) is used in the Markovian setting, we assume that this series converges on the class of functions under consideration.  All spectral
statements concerning \(B_k\) are understood on this domain. Although the Green--Kubo formula is often studied through Poisson equations and martingale approximations, our use of it here is as a covariance operator on finite-block observables.


\subsection{The i.i.d.\ nilpotent framework: a brief recap}

In the special case where the underlying sequence $\{X_t\}$ is i.i.d., the
lifted $k$-block chain possesses a remarkable finite-range structure
developed in Alhakim \cite{Alhakim2004} and Alhakim, Kawczak, and Molchanov \cite{AlhakimKawczakMolchanov2004}. After centering, the
window-chain operator becomes nilpotent in the sense that
\[
(P_k-\Pi)^k=0,
\]
where $\Pi$ denotes projection onto constants. As a consequence, on the subspace of centered functions the
Green--Kubo series (\ref{E:GreenKubo}) truncates after finitely many terms, and the asymptotic
covariance operator reduces to the finite polynomial form
\[
B_k
=
I+P_k+P_k^\ast+\cdots+P_k^{k-1}+(P_k^\ast)^{k-1}.
\]
This finite-range structure forces the covariance operator to possess the
integer spectrum
\[
\{0,1,\ldots,k\},
\]
leading to an orthogonal decomposition of $L^2(\pi_k)$ into overlap-induced
eigenspaces.

This integer spectrum underlies the independent-case decomposition developed in
\cite{AA2026}.  There, a general function on the \(k\)-window state space is decomposed into
orthogonal overlap-induced components, and the eigenvalue-one component is identified as the
incremental part, namely the information introduced at block length \(k\) beyond what
is representable using smaller block lengths. More generally, span-\(\ell\) block components correspond to the eigenvalue
\(k-\ell+1\), while differences between shifted copies of the same block component
belong to the kernel.

This independent decomposition is related in spirit to Hoeffding-type decompositions, in that both use orthogonality and successive averaging. The organizing principle, however, is different. In the overlap setting, the successive averaging is performed over side variables of consecutive blocks, so that components are indexed by block span and by their behavior under shifts, rather than by interaction order among unordered coordinates.


\section{Main results}

The preceding section recalled the Green--Kubo covariance operator for lifted \(k\)-block 
chains and the nilpotent overlap structure of the independent case. We now state 
the structural results that remain valid, or can be explicitly described, when 
the underlying sequence is Markovian. Proofs are deferred to Section~4.

We begin with the intrinsic incremental component. This component is defined before any use of reversibility, spectral calculus, or even the Markov property. Its definition is purely Hilbert-space theoretic. Indeed, it is the part of a \(k\)-block observable that is orthogonal to all information already present in the two adjacent \((k-1)\)-blocks.


Let $Y_t^{(k)}$ be the lifted \(k\)-block chain at step $t$. Define the adjacent
left and right \((k-1)\)-block spaces by
\[
W_{k-1}^L
=
\left\{
u(X_1,\ldots,X_{k-1}) :
u\in L^2_0(\pi_{k-1}^L)
\right\},
\]
and
\[
W_{k-1}^R
=
\left\{
v(X_2,\ldots,X_k) :
v\in L^2_0(\pi_{k-1}^R)
\right\},
\]
viewed as subspaces of \(\mathcal H_k\). Here \(\pi_{k-1}^L\) and \(\pi_{k-1}^R\) denote
the laws of the left and right \((k-1)\)-blocks, respectively. By stationarity
these laws coincide, but the two spaces are embedded differently inside \(\mathcal H_k\).
We set
\[
W_{k-1}
=
\overline{W_{k-1}^L+W_{k-1}^R},
\]
where the closure is taken in \(\mathcal H_k\). The closure is included so that
\(W_{k-1}\) is a closed subspace of \(\mathcal H_k\), allowing the orthogonal decomposition $\mathcal H_k=W_{k-1}\oplus W_{k-1}^{\perp}$.

\begin{definition}[Intrinsic incremental space]\label{D:orthogDef}
The incremental space at block length \(k\) is
\[
\mathcal I_k
=
W_{k-1}^{\perp}
=
\left(W_{k-1}^L+W_{k-1}^R\right)^\perp .
\]
Thus \(\mathcal I_k\) consists of those centered \(k\)-block observables that
are orthogonal to every observable based on either adjacent \((k-1)\)-block.
\end{definition}

Equivalently, by the defining property of conditional expectation, for $f\in\mathcal H_k$,
\begin{equation}\label{E:CondExpDefinition}
\mathcal I_k
=
\left\{
f :
E(f\mid X_1,\ldots,X_{k-1})=0,\ 
E(f\mid X_2,\ldots,X_k)=0
\right\}.
\end{equation}
\begin{remark}[An \(L^1\) version of the incremental condition]
The conditional-expectation characterization of the incremental space does not
require square integrability.  For \(f\in L^1_0(\pi_k)\), one may still define $\mathcal I_k$  by (\ref{E:CondExpDefinition}).
For a stationary Markov chain the same argument gives Equality (\ref{E:intersecKernels}) below
as an identity in \(L^1\). In the present paper we work in \(L^2\), since the
Green--Kubo covariance operator, spectral decomposition, and variance formulas
require the Hilbert-space structure.
\end{remark}
This formulation shows that the definition of \(\mathcal I_k\) is intrinsic and
does not rely on independence, reversibility, or the Markov property. In the current $L^2$ context it is the
orthogonal complement of all information already visible from the two adjacent
\((k-1)\)-blocks. 

When the underlying process is Markov, this static Hilbert-space definition acquires
a dynamic interpretation through the transition operator of the lifted \(k\)-block
chain. With \(B_k\) denoting the Green--Kubo covariance operator defined in
(\ref{E:GreenKubo}), we have the following result.

\begin{theorem}[Universal survival of the incremental component]\label{T:universal}
Let \(\{X_t\}_{t\in\mathbb Z}\) be a stationary Markov chain. Then the intrinsic incremental space satisfies 
\begin{equation}\label{E:intersecKernels}
\mathcal I_k=\ker P_k\cap \ker P_k^* .
\end{equation}
Consequently, the Green--Kubo series converges trivially on \(\mathcal I_k\), and,
\[
B_k f=f,\qquad f\in \mathcal I_k .
\]
Hence 1 is an eigenvalue of $B_k$ whenever \(\mathcal I_k\neq\{0\}\), and every element of
\(\mathcal I_k\) is an eigenfunction associated with this eigenvalue.
\end{theorem}

Using the language of prediction, the identity (\ref{E:intersecKernels}) provides a useful
interpretation of the incremental space.  The condition \(P_k f=0\) means that
the observable \(f\) has zero conditional expectation given the next window
state, while \(P_k^* f=0\) gives the corresponding condition for the time-reversed
chain.  Thus \(\mathcal I_k\) may be viewed as a two-sided innovation space.  In the present
overlapping-block setting, this abstract interpretation takes a concrete
statistical form: an observable \(f\in\mathcal I_k\) is invisible, in conditional mean,
from either adjacent \((k-1)\)-block.

Equivalently, \(\mathcal I_k\) measures the obstruction to reducing \(k\)-block observables
to adjacent \((k-1)\)-block information.  Indeed, by definition~\ref{D:orthogDef}, \(\mathcal I_k=\{0\}\) if and only if
\[
        \mathcal H_k=\overline{W^L_{k-1}+W^R_{k-1}}.
\]
That is, every centered \(k\)-block observable can be approximated in \(L^2\) by
sums of observables depending on the left and right adjacent \((k-1)\)-blocks.
Thus a nonzero \(\mathcal I_k\) records precisely the presence of directions in \(\mathcal H_k\)
that are not reducible to adjacent shorter blocks.  In the Markov case,
Theorem~\ref{T:universal} identifies this obstruction dynamically.

In finite-alphabet or symbolic settings, this provides the operator-theoretic
counterpart to the order-selection and context-tree viewpoints discussed in the
introduction: here the reduction question is posed at the level of block
observables and covariance directions rather than conditional probability models.

We next examine the first nontrivial block length, \(k=2\), where this abstract
definition admits a complete and useful decomposition.  In this case a two-point
observable separates into a part explained by one-coordinate marginal information
and a pure transition component belonging to \(\mathcal I_2\).  This two-point separation
will serve as the base case for the later decomposition of higher-block
observables.

\begin{theorem}(The case \(k=2\): separation of marginal and transition information).\label{T:Casek=2}
Let \(\{X_t\}_{t\in\mathbb Z}\) be a stationary reversible Markov chain on a state
space \((E,\mathcal E)\), with invariant law \(\nu\) and transition
kernel \(K\). 
Assume that \(I-K^2\) is invertible on \(L^2_0(\nu)\), and let \(f\in \mathcal H_2\).
Then there exist unique functions \(a,b\in L^2_0(\nu)\), explicitly computable, and a unique function \(g\in\mathcal I_2\) such
that
\[
f(x,y)=g(x,y)+a(x)+b(y).
\]
Moreover, $g(x,y)$ and $a(x)+b(y)$ are orthogonal.

\end{theorem}
Thus \(a(x)+b(y)\) is the one-coordinate marginal part of \(f\), while
\(g(x,y)\) is the intrinsic incremental, or pure transition, component.
Note that for the special case of a finite-state ergodic chain, the assumption that \(I-K^2\) is invertible amounts to the condition that the
spectrum of \(K\) on \(L^2_0(\nu)\) does not include \(-1\).
This \(k=2\) decomposition separates the first genuinely two-point contribution from the one-coordinate marginal layer. 
For larger \(k\), the same principle appears recursively: lower-span incremental
components may be embedded into longer windows.

\begin{theorem}(Shift-inherited eigenvalues).\label{T:Shift-inherited}
Let \(\{X_t\}_{t\in\mathbb Z}\) be a stationary Markov chain. Let \(2\le m\le k\) and let $h\in \mathcal I_m$.
Define
\[
\widehat{h}_{m,k}(x_1,\ldots,x_k)
=
\sum_{i=1}^{k-m+1}
h(x_i,\ldots,x_{i+m-1}).
\]
For this function the Green--Kubo series is finite, and
\[
B_k \widehat{h}_{m,k}=(k-m+1) \widehat{h}_{m,k}.
\]
\end{theorem}

Taking \(m=k,k-1,\ldots,2\), the theorem gives the eigenvalues $1,2,\ldots,k-1$.
We refer to these as the shift-inherited eigenvalues.

\begin{remark}[Block sums vs. translation contrasts]
The theorem identifies the informative block-sum direction generated by the valid
translates of \(h\in\mathcal I_m\).  It should not be read as saying that the whole
linear span of these translates is an eigenspace for the eigenvalue \(k-m+1\).
Indeed, the raw translate span also contains covariance-kernel directions.  For
example, when \(m=k-1\), the two translates
\[
h(x_1,\ldots,x_{k-1}),\qquad h(x_2,\ldots,x_k)
\]
produce the eigenvalue-two direction
\[
h(x_1,\ldots,x_{k-1})+h(x_2,\ldots,x_k),
\]
whereas the contrast
\[
h(x_1,\ldots,x_{k-1})-h(x_2,\ldots,x_k)
\]
is a kernel, or coboundary, direction.  More generally, the integer eigenvalue
comes from the block-sum direction, while translation contrasts belong to the
zero-variance part.
\end{remark}

The preceding theorem accounts for the integer eigenvalues
\(1,\ldots,k-1\), and does so without assuming reversibility. The remaining layer is different: it is
the one-coordinate marginal layer, where the covariance structure depends on
the transition dynamics of the underlying Markov chain.  Under reversibility,
this layer can be diagonalized through the spectral resolution of the base
transition operator.

\begin{theorem}(Reversible marginal spectrum).\label{T:marginal}\\
Assume that the stationary Markov chain is reversible, and let \(K\) denote the
transition operator of the base chain acting on \(L^2_0(\nu)\).  For every 
\(\lambda\in Spec(K)\backslash\{1\}\), define
\begin{equation}\label{marginalEigenvalues}
\theta_k(\lambda)
=
\frac{k-(k-2)\lambda}{1-\lambda}.
\end{equation}
Then, on the marginal subspace generated by one-coordinate functions in the  \(k\)-window space, the operator 
\(B_k\)  admits \(\theta_k(\lambda)\) as eigenvalues. More precisely, 
 if \(h\in L^2_0(\nu)\) satisfies \(Kh=\lambda h\), then the boundary-corrected function
\[
M_{\lambda,k}(x_1,\ldots,x_k)
=
\frac{1}{1-\lambda}h(x_1)
+
\sum_{i=2}^{k-1}h(x_i)
+
\frac{1}{1-\lambda}h(x_k)
\]
is an eigenfunction of \(B_k\) with eigenvalue \(\theta_k(\lambda)\).
\end{theorem}


\begin{remark}
When the base sequence is independent, the transition operator satisfies
\(K=0\) on \(L^2_0(\nu)\).  Hence
\[
\theta_k(0)=k,
\]
and the transition-inherited marginal layer reduces to the final integer
eigenvalue \(k\) of the classical nilpotent spectrum $\{0,1,\ldots,k\}$. In the current case,  the family \(\theta_k(\lambda)\) records the Markovian dynamics of the base chain.
\end{remark}


\subsection{Downward decomposition}

We now state the structural decomposition of a general window observable.  The
one-coordinate layer is determined first, while the remaining nonzero layers are
generated recursively by incremental components of smaller spans.

\begin{proposition}(Projection onto the marginal subspace).\label{P:projectionMarginal}\\
Assume reversibility and that \(I-K^2\) is invertible on
\(L^2_0(\nu)\).  Let
\[
\mathcal M_k
=
\overline{
\left\{
a_1(X_1)+\cdots+a_k(X_k):
a_i\in L^2_0(\nu)
\right\}}
\subset \mathcal H_k
\]
be the one-variable marginal space.  For \(f\in \mathcal H_k\), define
\[
q_i(x)=E\{f(X_1,\ldots,X_k)\mid X_i=x\},
\qquad i=1,\ldots,k,
\]
and write
\[
q=(q_1,\ldots,q_k)^T.
\] 
Then the orthogonal projection of \(f\) onto \(\mathcal M_k\) is the unique element 
\(a_1(X_1)+\cdots+a_k(X_k)\), 
whose coefficient vector \((a_1,\ldots,a_k)^T\) solves the system
\[
(a_1,\ldots,a_k)^T=
(I-K^2)^{-1}
\begin{pmatrix}
I & -K & 0 & \cdots & 0\\
-K & I+K^2 & -K & \ddots & \vdots\\
0 & -K & I+K^2 & \ddots & 0\\
\vdots & \ddots & \ddots & \ddots & -K\\
0 & \cdots & 0 & -K & I
\end{pmatrix}
q ,
\]
where the factor \((I-K^2)^{-1}\) is applied componentwise.
\end{proposition}

For spans at least three, the projection onto the incremental space has
a simple inclusion--exclusion form.

\begin{proposition}[Intrinsic projection for spans \(m\ge 3\)]\label{P:IntrinsicProjection}
Let \(m\ge 3\), and let $\psi\in\mathcal H_m= L^2_0(\pi_m)$. Define
\[
L_m\psi
=
E(\psi\mid X_1,\ldots,X_{m-1}),
\qquad
R_m\psi
=
E(\psi\mid X_2,\ldots,X_m),
\]
and
\[
C_m\psi
=
E(\psi\mid X_2,\ldots,X_{m-1}).
\]
Then the orthogonal projection of \(\psi\) onto \(\mathcal I_m\) is
\[
\operatorname{Proj}_{\mathcal I_m}\psi
=
\psi-L_m\psi-R_m\psi+C_m\psi .
\]
\end{proposition}

The case \(m=2\) is the two-point separation theorem stated above; it replaces
the middle-block correction by the one-coordinate marginal projection.

The preceding proposition is local. It extracts the intrinsic incremental part of a single \(m\)-block observable. The downward decomposition below applies this local purification successively across the possible spans \(m=k,k-1,\ldots,2\), embedding each purified component back into the original \(k\)-block space. Let \(\mathcal M_k\) be the one-coordinate marginal space introduced in Proposition~\ref{P:projectionMarginal}. We write \[ \mathcal M_k = \mathcal M_k^+\oplus(\mathcal M_k\cap\ker B_k), \] where \(\mathcal M_k^+\) denotes the nonzero spectral part of \(\mathcal M_k\). In the reversible case, \(\mathcal M_k^+\) is the marginal spectral layer diagonalized through the boundary-corrected eigenfunctions described in Theorem~\ref{T:marginal}.

For \(2\le m\le k\), define \[ S_{m,k} = \left\{ \sum_{i=1}^{k-m+1} h(x_i,\ldots,x_{i+m-1}) : h\in \mathcal I_m \right\}. \]

\begin{theorem}(Downward spectral decomposition).\label{T:downwardDecomp}\\
Assume that the chain is reversible, that \(I-K^2\) is invertible on
\(L^2_0(\nu)\), and that the Green--Kubo series defining \(B_k\) converges on all of \(H_k\).  Then
\[
        \mathcal H_k
        =
        \mathcal M_k^+
        \oplus S_{2,k}\oplus\cdots\oplus S_{k,k}\oplus\ker B_k .
\]
Consequently, every \(f\in \mathcal H_k\) has a unique decomposition
\begin{equation}\label{E:decomp}
        f=f_{\rm marg}+\sum_{m=2}^k f^{[m]}+f_0,
\end{equation}
where
\[
        f_{\rm marg}\in \mathcal M_k^+,\qquad
        f^{[m]}\in S_{m,k},\qquad
        f_0\in\ker B_k .
\]
\end{theorem}


\begin{corollary}[Variance decomposition]
Under the hypotheses of Theorem~\ref{T:downwardDecomp}, let $f$ be the decomposition (\ref{E:decomp}).  Then
\[
\sigma_k^2(f)
=
\langle B_k f,f\rangle
=
\sigma^2_{\mathrm{marg}}(f)
+
\sum_{m=2}^k \sigma^2_{[m]}(f),
\]
where
\[
\sigma^2_{[m]}(f)
=
\langle B_k f^{[m]},f^{[m]}\rangle
=
(k-m+1)\|f^{[m]}\|^2 =(k-m+1)^2\|h_m\|^2 .
\]
In particular,
\[
\sigma^2_{\mathrm{inc}}(f)
=
\sigma^2_{[k]}(f)
=
\|f^{[k]}\|^2 .
\]
\end{corollary}
We use the notation
\[
f^{[m]}=\sum_{i=1}^{k-m+1}h_m(x_i,\ldots,x_{i+m-1}),
\qquad h_m\in \mathcal I_m .
\]
This corollary follows immediately from Theorem~\ref{T:downwardDecomp}, so we do
not give it a separate proof. We only note that the rightmost equality for
\(\sigma^2_{[m]}(f)\) uses the fact that the shifted copies of \(h_m\) are pairwise
orthogonal; this follows from the Markov property and the two-sided innovation
conditions defining \(\mathcal I_m\).


Under the spectral resolution of the reversible base operator \(K\), the marginal contribution \(\sigma^2_{\mathrm{marg}}(f)%
=\langle B_k f_{\mathrm{marg}},f_{\mathrm{marg}}\rangle\), 
is represented as
\[
\sigma^2_{\mathrm{marg}}(f)
=
\int_{\operatorname{Spec}(K)\backslash\{1\}}
\theta_k(\lambda)\, d\mu_f(\lambda),
\]
where \(\theta_k(\lambda)\) is given in (\ref{marginalEigenvalues})
and \(\mu_f\) is the spectral measure of the marginal component under the
marginal-layer identification with \(K\).


\section{Proofs of the main results}

\subsection{Proof of Theorem \ref{T:universal}}

We prove first the identity
\[
\mathcal I_k=\ker P_k\cap\ker P_k^* .
\]
By stationarity, it is enough to work with $Y_1^{(k)}=(X_1,\ldots,X_k)$.
Since \(P_k\) is the transition operator of the lifted chain,
\[
P_k f(Y_1^{(k)})
=
E\{f(Y_2^{(k)})\mid Y_1^{(k)}\}.
\]
Expanding $Y_1^{(k)}$ and $Y_2^{(k)}$ we obtain
\[
P_k f(X_1,\ldots,X_k)
=
E\{f(X_2,\ldots,X_{k+1})\mid X_1,\ldots,X_k\}.
\]
By the Markov property, the right-hand side depends only on
\((X_2,\ldots,X_k)\). Hence \(P_k f=0\) is equivalent, by stationarity, to
\[
E\{f(X_1,\ldots,X_k)\mid X_1,\ldots,X_{k-1}\}=0.
\]
Similarly,
\[
P_k^* f(Y_1^{(k)})
=
E\{f(Y_0^{(k)})\mid Y_1^{(k)}\}.
\]
That is,
\[
P_k^* f(X_1,\ldots,X_k)
=
E\{f(X_0,\ldots,X_{k-1})\mid X_1,\ldots,X_k\}.
\]
By the Markov property for the reversed stationary chain, this depends only on
\((X_1,\ldots,X_{k-1})\). Hence \(P_k^*f=0\) is equivalent, by stationarity, to
\[
E\{f(X_1,\ldots,X_k)\mid X_2,\ldots,X_k\}=0.
\]
Combining the two equivalences, we get
\[
\ker P_k\cap\ker P_k^*
=
\left\{
f\in \mathcal H_k:
E(f\mid X_1,\ldots,X_{k-1})=0,\ 
E(f\mid X_2,\ldots,X_k)=0
\right\}.
\]
By the conditional-expectation characterization of \(\mathcal I_k\). The right-hand
side is exactly \(\mathcal I_k\), by the conditional expectation characteriztion in Definition 3.1.  This establishes Equation~(\ref{E:intersecKernels}).

Finally, if \(f\in\mathcal I_k\), then \(P_kf=0\) and \(P_k^*f=0\).  Therefore
\[
P_k^h f=0,\qquad (P_k^*)^h f=0,\qquad h\ge1.
\]
Consequently, $\mathcal I_k$ is a subset of the domain of definition of  the Green--Kubo operator, and 
\[
B_kf
=
\left(
I+\sum_{h\ge1}P_k^h+\sum_{h\ge1}(P_k^*)^h
\right)f
=
f.
\]
This completes the proof.

\subsection{Proof of Theorem \ref{T:Casek=2}}

Let
\[
L f(x)
=
E\{f(X_1,X_2)\mid X_1=x\},
\qquad
R f(y)
=
E\{f(X_1,X_2)\mid X_2=y\}.
\]
Since \(f\in L^2_0(\pi_2)\), both \(Lf\) and \(Rf\) are centered elements of
\(L^2_0(\nu)\).  We seek centered functions \(a,b\in L^2_0(\nu)\) such that
\begin{equation}\label{E:I2Comp}
g(x,y)=f(x,y)-a(x)-b(y)
\end{equation}
belongs to \(\mathcal I_2\).  By the conditional-expectation characterization of
\(\mathcal I_2\), this is equivalent to
\begin{equation}\label{E:Con_k=2}
E(g(X_1,X_2)\mid X_1)=0,
\qquad
E(g(X_1,X_2)\mid X_2)=0.
\end{equation}
The first condition,  Equation~(\ref{E:I2Comp}), and the operator representation (\ref{E:MarkovOperatorK}) give
\[
a+Kb=Lf.
\]
Since the chain is reversible,
\[
E\{a(X_1)\mid X_2=y\}=(Ka)(y).
\]
Thus the second condition of (\ref{E:Con_k=2}) gives
\[
Ka+b=Rf.
\]

Thus the marginal part \(a(x)+b(y)\) is determined by the system
\[
\begin{cases}
a+Kb=Lf,\\
Ka+b=Rf.
\end{cases}
\]
We solve this system by elimination.  From the first equation,
\[
        a=Lf-Kb.
\]
Substituting into the second equation gives
\[
        K(Lf-Kb)+b=Rf,
\]
and hence
\[
        (I-K^2)b=Rf-KLf.
\]
Similarly, eliminating \(b\) gives
\[
        (I-K^2)a=Lf-KRf.
\]
Under the stated nondegeneracy assumption that \(I-K^2\) is invertible on
\(L^2_0(\nu)\), these equations determine \(a\) and \(b\) uniquely,
\begin{equation}\label{E:a+b}
\begin{cases}
        a=(I-K^2)^{-1}(Lf-KRf),\\
        b=(I-K^2)^{-1}(Rf-KLf).
\end{cases}
\end{equation}
Thus the marginal part \(a(x)+b(y)\) is uniquely determined, and the residual
\(g=f-a-b\) belongs to \(\mathcal I_2\).

It remains to prove uniqueness.  Suppose that
\[
f=g+a(X_1)+b(X_2)
=
\tilde g+\tilde a(X_1)+\tilde b(X_2)
\]
are two such decompositions.  Set $\alpha=a-\tilde a$, and $\beta=b-\tilde b$ so that
\[
(g-\tilde g)+\alpha(X_1)+\beta(X_2)=0.
\]
The first term lies in \(\mathcal I_2\), while the second lies in
\(W_1^L+W_1^R\).  Since these spaces are orthogonal, \(g=\tilde g\), and hence
\[
\alpha(X_1)+\beta(X_2)=0.
\]
Taking conditional expectations with respect to \(X_1\) and \(X_2\) gives
\[
\alpha+K\beta=0,
\qquad
K\alpha+\beta=0.
\]
Therefore
\[
(I-K^2)\alpha=0,
\qquad
(I-K^2)\beta=0.
\]
By the assumed invertibility of \(I-K^2\) on \(L^2_0(\nu)\), we obtain
\[
\alpha=\beta=0.
\]


\begin{remark}
This theorem shows that, already at block length $k=2$, one can separate the
information carried by a two-point window into:
\[
\text{marginal information } a(x)+b(y)
\quad\text{and}\quad
\text{pure transition information } g(x,y).
\]

This is formally analogous to the first Hoeffding decomposition of a two-variable
kernel.  For dependent observations, such Hoeffding decompositions have been used
effectively in the analysis of \(U\)-statistics; see, for example, Dehling and Wendler
\cite{DehlingWendler2010}.  The distinction here is that the residual \(g\) is defined
relative to the Markov transition law \(\nu(dx)K(x,dy)\), rather than relative to the
product marginal law used in the classical Hoeffding projection. Thus \(g\) is indeed 
a transition-level incremental component, and serves as the basic model for higher window sizes.
\end{remark}

\subsection{Proof of Theorem \ref{T:Shift-inherited}}

Let
\[
\ell=k-m+1,
\]
and, for \(i=1,\ldots,\ell\), define the \(i\)-th valid translate of \(h\) inside a
\(k\)-block by
\[
e_i(h)(x_1,\ldots,x_k)
=
h(x_i,\ldots,x_{i+m-1}).
\]
Thus
\[
\widehat{h}_{m,k}=\sum_{i=1}^{\ell} e_i(h).
\]

We first record the action of the operator \(P_k\) on these
translates. By Equation (\ref{E:MarkovOperatorP}) it is immediate that
\[
P_k e_i(h)=e_{i+1}(h),\qquad i=1,\ldots,\ell-1.
\]
The last translate includes the variable \(x_k\), so shifting the block introduces
a new variable drawn according to the transition kernel \(K(x_k,dz)\).  Explicitly,
\[
(P_k e_\ell(h))(x_1,\ldots,x_k)
=
\int_E h(x_{\ell+1},\ldots,x_k,z)\,K(x_k,dz).
\]
With the lowercase variables understood as fixed, this is equivalent to
\[
(P_k e_\ell(h))(x_1,\ldots,x_k)
=
E\{h(x_{\ell+1},\ldots,x_k,X_{k+1})\mid X_k=x_k\}.
\]

By the Markov property and the defining conditional-mean-zero property of
\(\mathcal I_m\), this expression vanishes. Hence \(P_k e_\ell(h)=0\).

Therefore
\begin{equation}\label{E:Pforward}
P_k e_i(h)=
\begin{cases}
e_{i+1}(h), & 1\le i\le \ell-1,\\
0, & i=\ell.
\end{cases}
\end{equation}

The adjoint calculation is analogous, using the reversed lifted chain. As a result,
\begin{equation}\label{E:Pbackward}
P_k^* e_i(h)=
\begin{cases}
0, & i=1,\\
e_{i-1}(h), & 2\le i\le \ell.
\end{cases}
\end{equation}

It follows that, for \(r=1,\ldots,\ell-1\),
\[
P_k^r \widehat{h}_{m,k}
=
e_{r+1}(h)+e_{r+2}(h)+\cdots+e_\ell(h),
\]
and
\[
(P_k^*)^r \widehat{h}_{m,k}
=
e_1(h)+e_2(h)+\cdots+e_{\ell-r}(h).
\]
For \(r\ge \ell\), both terms vanish.
Therefore, using the Green--Kubo covariance operator,
\[
B_k\widehat{h}_{m,k}
=
\widehat{h}_{m,k}
+
\sum_{r=1}^{\ell-1} P_k^r\,\widehat{h}_{m,k}
+
\sum_{r=1}^{\ell-1} (P_k^*)^r\,\widehat{h}_{m,k}.
\]
We now collect the coefficient of a fixed translate \(e_j(h)\).  It appears once in
\(\widehat{h}_{m,k}\), appears in \(P_k^r\widehat{h}_{m,k}\) for
\[
r=1,\ldots,j-1,
\]
and appears in \((P_k^*)^r\widehat{h}_{m,k}\) for
\[
r=1,\ldots,\ell-j.
\]
Hence the total coefficient of \(e_j(h)\) is
\[
1+(j-1)+(\ell-j)=\ell.
\]
Since this holds for every \(j=1,\ldots,\ell\), we obtain
\[
B_k\widehat{h}_{m,k}=\ell \widehat{h}_{m,k}.
\]
Thus the embedded block-sum generated by \(h\in\mathcal I_m\) is an eigenfunction
of \(B_k\) with eigenvalue \(\ell=k-m+1\), as claimed.

\subsection{Proof of Theorem \ref{T:marginal}}


Let \(h\in L^2_0(\nu)\) be an eigenfunction of the base transition operator \(K\), say
\begin{equation}\label{E:eigenfunctionK}
Kh=\lambda h,\qquad \lambda\neq 1.
\end{equation}
Note that $h$ is a one-variable function and write 
\begin{equation}\label{E:SlidingOne}
e_i=e_i(h)(x_1,\ldots,x_k)=h(x_i),
\qquad 
i=1,\ldots,k.
\end{equation}
By the block-shift formula established in the proof of the previous theorem,
\[
P_k e_i=e_{i+1},\qquad i=1,\ldots,k-1.
\]
At the right boundary, the shift introduces a new variable drawn from \(K(x_k,dz)\).
Hence, (\ref{E:eigenfunctionK}) and  (\ref{E:SlidingOne}) and this observation imply
\[
P_k e_k=\lambda e_k.
\]
Similarly, by reversibility,
\[
P_k^*e_i=e_{i-1},\qquad i=2,\ldots,k,
\qquad
P_k^*e_1=\lambda e_1.
\]

It follows iteratively that
\[
P_k^r e_i
=
\begin{cases}
e_{i+r}, & i+r\le k,\\[4pt]
\lambda^{\,i+r-k}e_k, & i+r>k,
\end{cases}
\qquad
(P_k^*)^r e_i
=
\begin{cases}
e_{i-r}, & i-r\ge 1,\\[4pt]
\lambda^{\,r-i+1}e_1, & i-r<1.
\end{cases}
\qquad r\ge1.
\]

Therefore,
\[
\begin{aligned}
B_k e_i
&=
e_i+\sum_{r\ge1}P_k^r e_i+\sum_{r\ge1}(P_k^*)^r e_i  \\
&=
\sum_{j=1}^k e_j
+
\left(\sum_{s\ge1}\lambda^s\right)(e_1+e_k) \\
&=
\sum_{j=1}^k e_j
+
\frac{\lambda}{1-\lambda}(e_1+e_k) \\
&=
\frac{1}{1-\lambda}e_1
+
\sum_{j=2}^{k-1}e_j
+
\frac{1}{1-\lambda}e_k.
\end{aligned}
\]
The right-hand side is independent of \(i\). We will denote it by $M_{\lambda,k}$. Thus, referring to Equation~(\ref{E:SlidingOne}) we now have, for $i=1,\ldots,k$
\[
B_k e_i
=
M_{\lambda,k}(x_1,\ldots,x_k)
=
\frac{1}{1-\lambda}h(x_1)
+
\sum_{j=2}^{k-1}h(x_j)
+
\frac{1}{1-\lambda}h(x_k).
\]

It follows that for any marginal combination $v=\displaystyle\sum_{i=1}^k c_i e_i$ we have
\[
B_k v
=
\left(\sum_{i=1}^k c_i\right)M_{\lambda,k}.
\]
In particular, all marginal contrasts satisfying $\displaystyle\sum_{i=1}^k c_i=0$ 
belong to \(\ker B_k\).  On the other hand, for \(v=M_{\lambda,k}\), the sum of the
coefficients is
\[
(k-2)+\frac{2}{1-\lambda}
=
\frac{k-(k-2)\lambda}{1-\lambda}.
\]
Therefore
\(M_{\lambda,k}\) is an eigenfunction of \(B_k\) with the claimed eigenvalue $\theta_k(\lambda)$.



\subsection{Proof of Proposition \ref{P:projectionMarginal}}

Let $f\in\mathcal H_k$ and write
the orthogonal projection of \(f\) onto \(\mathcal M_k\) as
\[
f_{\mathrm{marg}}
=
a_1(X_1)+\cdots+a_k(X_k).
\]
This projection is characterized by the normal equations
\[
f-f_{\mathrm{marg}}\perp \mathcal M_k.
\]
Equivalently, $E\{f-f_{\mathrm{marg}}\mid X_i\}=0$ for each \(i=1,\ldots,k\).
By the definition of \(q_i\), this gives
\[
q_i
=
E\{f\mid X_i\}
=
\sum_{j=1}^k E\{a_j(X_j)\mid X_i\}.
\]
This is the same projection principle that underlies finite linear prediction:
the best additive marginal approximation is characterized by orthogonality of
the residual to each of the one-coordinate subspaces. In the present Markov
setting, the resulting normal equations are coupled through the transition
powers of \(K\).
Since the chain is stationary and reversible, the conditional expectation of
\(a_j(X_j)\) given \(X_i\) is
\[
E\{a_j(X_j)\mid X_i\}
=
K^{|i-j|}a_j.
\]
which we obtain by iterating (\ref{E:MarkovOperatorK}) $|i-j|$ times. 
 Hence the
normal equations become
\[
\sum_{j=1}^k K^{|i-j|}a_j=q_i,
\qquad i=1,\ldots,k.
\]
In vector form,
\[
\begin{pmatrix}
I & K & K^2 & \cdots & K^{k-1}\\
K & I & K & \cdots & K^{k-2}\\
K^2 & K & I & \cdots & K^{k-3}\\
\vdots & \vdots & \vdots & \ddots & \vdots\\
K^{k-1} & K^{k-2} & K^{k-3} & \cdots & I
\end{pmatrix}
\begin{pmatrix}
a_1\\
a_2\\
\vdots\\
a_k
\end{pmatrix}
=
\begin{pmatrix}
q_1\\
q_2\\
\vdots\\
q_k
\end{pmatrix}.
\]
The operator matrix on the left is the operator-valued analogue of the classical Kac--Murdock--Szegő
Toeplitz matrix \((\rho^{|i-j|})\), whose inverse is tridiagonal; see
\cite{KacMurdockSzego1953}.  Since the entries here are all powers of the same
operator \(K\), all entries commute and the same algebraic identity applies, with \(\rho\) replaced by \(K\).
Under the assumed invertibility of \(I-K^2\), its inverse is precisely the
tridiagonal operator matrix displayed in the statement of the proposition.
It follows that the functions \(a_1,\ldots,a_k\) are uniquely determined as stated.


\subsection{Proof of Proposition \ref{P:IntrinsicProjection}}

Define the $\sigma$-fields
\[
\mathcal F_L=\sigma(X_1,\ldots,X_{m-1}),\qquad
\mathcal F_R=\sigma(X_2,\ldots,X_m),
\]
and
\[
\mathcal F_C=\sigma(X_2,\ldots,X_{m-1}).
\]
Thus
\[
L_m\psi=E(\psi\mid \mathcal F_L),\qquad
R_m\psi=E(\psi\mid \mathcal F_R),\qquad
C_m\psi=E(\psi\mid \mathcal F_C).
\]

The Markov property implies that the two endpoint variables \(X_1\) and \(X_m\)
are conditionally independent given the middle block
\[
(X_2,\ldots,X_{m-1}).
\]
Equivalently, the two sigma-fields \(\mathcal F_L\) and \(\mathcal F_R\) are
conditionally independent given their common sub-sigma-field \(\mathcal F_C\).
Consequently,
\[
E\{R_m\psi\mid \mathcal F_L\}
=
E\{E(\psi\mid \mathcal F_R)\mid \mathcal F_L\}
=
E(\psi\mid \mathcal F_C)
=
C_m\psi,
\]
and similarly
\[
E\{L_m\psi\mid \mathcal F_R\}=C_m\psi.
\]

Now define
\[
\psi^{\mathrm{inc}}
=
\psi-L_m\psi-R_m\psi+C_m\psi.
\]
We first show that \(\psi^{\mathrm{inc}}\in\mathcal I_m\).  Conditioning on the left
\((m-1)\)-block gives
\[
\begin{aligned}
E(\psi^{\mathrm{inc}}\mid \mathcal F_L)
&=
E(\psi\mid \mathcal F_L)
-
E(L_m\psi\mid \mathcal F_L)
-
E(R_m\psi\mid \mathcal F_L)
+
E(C_m\psi\mid \mathcal F_L)  \\
&=
L_m\psi-L_m\psi-C_m\psi+C_m\psi\\
&=0.
\end{aligned}
\]
Similarly, conditioning on the right \((m-1)\)-block gives
\[
\begin{aligned}
E(\psi^{\mathrm{inc}}\mid \mathcal F_R)
&=
E(\psi\mid \mathcal F_R)
-
E(L_m\psi\mid \mathcal F_R)
-
E(R_m\psi\mid \mathcal F_R)
+
E(C_m\psi\mid \mathcal F_R)\\
&=
R_m\psi-C_m\psi-R_m\psi+C_m\psi\\
&=0.
\end{aligned}
\]

It remains to identify this element as the orthogonal projection of \(\psi\) onto
\(\mathcal I_m\).  We have
\[
\psi-\psi^{\mathrm{inc}}
=
L_m\psi+R_m\psi-C_m\psi.
\]
Here \(L_m\psi\) is a function of \((X_1,\ldots,X_{m-1})\), while
\(R_m\psi-C_m\psi\) is a function of \((X_2,\ldots,X_m)\).  Hence
\[
\psi-\psi^{\mathrm{inc}}\in W_{m-1}.
\]
Since $\mathcal I_m=W_{m-1}^{\perp}$, 
we have decomposed \(\psi\) as
\[
\psi=\psi^{\mathrm{inc}}+(\psi-\psi^{\mathrm{inc}}),
\]
with
\[
\psi^{\mathrm{inc}}\in\mathcal I_m,
\qquad
\psi-\psi^{\mathrm{inc}}\in W_{m-1}.
\]
Thus \(\psi^{\mathrm{inc}}\) is the orthogonal projection of \(\psi\) onto
\(\mathcal I_m\).
This proves the proposition.

\subsection{Proof of Theorem \ref{T:downwardDecomp}}

Since the chain is reversible, \(B_k\) is self-adjoint on \(\mathcal H_k\).
Hence spectral subspaces corresponding to distinct spectral values are orthogonal.

We first organize \(\mathcal H_k\) by consecutive block span.  For \(1\le m\le k\),
let \(V_{m,k}\) be the closed linear span in \(\mathcal H_k\) of all observables
depending on \(m\) consecutive coordinates:
\[
        V_{m,k}
        =
        \overline{\operatorname{span}}
        \left\{
        \psi(x_i,\ldots,x_{i+m-1}):
        \psi\in \mathcal H_m,\ 1\le i\le k-m+1
        \right\}.
\]
Thus
\[
        V_{1,k}=\mathcal M_k,\qquad V_{k,k}=\mathcal H_k.
\]
For each \(m\ge2\), the local decomposition of an \(m\)-block observable gives
\[
        \mathcal H_m=W_{m-1}\oplus \mathcal I_m.
\]
Embedding this identity into each possible position inside a \(k\)-block gives
\[
        V_{m,k}=V_{m-1,k}+\widetilde S_{m,k},
\]
where
\[
        \widetilde S_{m,k}
        =
        \overline{\operatorname{span}}
        \left\{
        h(x_i,\ldots,x_{i+m-1}):
        h\in \mathcal I_m,\ 1\le i\le k-m+1
        \right\}.
\]
Iterating over \(m=2,\ldots,k\), we obtain
\[
        \mathcal H_k
        =
        \mathcal M_k+\sum_{m=2}^k \widetilde S_{m,k}.
\]

We next identify the nonzero covariance part of each raw translate space
\(\widetilde S_{m,k}\).  Fix \(m\), and set $\ell=k-m+1$.
For \(h\in\mathcal I_m\) and $ i=1,\ldots,\ell$, let $e_i(h)$ be as defined in the proof of Theorem~\ref{T:Shift-inherited}.
Therefore by Equations (\ref{E:Pforward}) and (\ref{E:Pbackward})
\[
        B_k e_i(h)=\sum_{j=1}^{\ell} e_j(h),
        \qquad i=1,\ldots,\ell .
\]
It follows that, for arbitrary \(h_1,\ldots,h_\ell\in\mathcal I_m\),
\[
        B_k\left(\sum_{i=1}^{\ell}e_i(h_i)\right)
        =
        \sum_{j=1}^{\ell}e_j\left(\sum_{i=1}^{\ell}h_i\right).
\]

Thus the only nonzero covariance contribution of $\widetilde{S}_{m,k}$  lies in the block-sum space  $S_{m,k}$; the remaining translate contrasts lie in \(\ker B_k\).
Indeed, if
\[
        \bar h=\frac1\ell\sum_{i=1}^{\ell}h_i,
\]
then
\[
        \sum_{i=1}^{\ell}e_i(h_i)
        =
        \sum_{i=1}^{\ell}e_i(\bar h)
        +
        \sum_{i=1}^{\ell}e_i(h_i-\bar h),
\]
where $\displaystyle\sum_{i=1}^{\ell}e_i(\bar h)\in S_{m,k}$,
and since $\displaystyle\sum_{i=1}^{\ell}(h_i-\bar h)=0$, 
\[
        \sum_{i=1}^{\ell}e_i(h_i-\bar h)\in\ker B_k,
\]
By Theorem \ref{T:Shift-inherited}, every \(u\in S_{m,k}\) is an eigenfunction of $B_k$ with eigenvalue $\ell$. 
Thus every element of the raw translate space \(\widetilde S_{m,k}\) decomposes
into a block-sum eigencomponent in \(S_{m,k}\) and a covariance-kernel contrast.
Equivalently,
\[
        \widetilde S_{m,k}
        \subseteq
        S_{m,k}+\ker B_k .
\]
Since \(S_{m,k}\subseteq \widetilde S_{m,k}\) by definition, the only additional
directions in \(\widetilde S_{m,k}\) are translation contrasts lying in \(\ker B_k\).
Therefore the span filtration gives
\[
        \mathcal H_k
        =
        \mathcal M_k+\sum_{m=2}^k S_{m,k}+\ker B_k .
\]

It remains to split the marginal layer. By Proposition~3.8 and Theorem~3.6, the
one-coordinate marginal space decomposes as
\[
       \mathcal  M_k=\mathcal M_k^+\oplus (\mathcal M_k\cap\ker B_k),
\]
where \(\mathcal M_k^+\) is the nonzero marginal spectral part described through the
boundary-corrected eigenfunctions of Theorem~3.6. Absorbing
\(\mathcal M_k\cap\ker B_k\) into the kernel term yields
\[
       \mathcal H_k
        =
       \mathcal  M_k^+
        +
        \sum_{m=2}^k S_{m,k}
        +
        \ker B_k .
\]

We now verify that this sum is orthogonal. First, the spaces \(S_{m,k}\) and \(S_{n,k}\) correspond to distinct eigenvalues
whenever \(m\ne n\), and are therefore orthogonal because \(B_k\) is self-adjoint.
The marginal nonzero part \(\mathcal M_k^+\) has spectral values
\[
        \theta_k(\lambda)
        =
        \frac{k-(k-2)\lambda}{1-\lambda},
        \qquad
        \lambda\in\operatorname{Spec}(K),\quad \lambda\ne1 .
\]
These values do not coincide with the shift-inherited integers
\(1,\ldots,k-1\), under the standing assumption that
\(-1\notin\operatorname{Spec}(K)\). Indeed, solving
\[
        \theta_k(\lambda)=j,
        \qquad
        j\in\{1,\ldots,k-1\},
\]
gives no admissible solution for \(j\le k-2\), while \(j=k-1\) gives precisely
\(\lambda=-1\). Hence \(\mathcal M_k^+\) is orthogonal to each \(S_{m,k}\). Finally,
\(\ker B_k\) is orthogonal to all nonzero spectral subspaces of the self-adjoint
operator \(B_k\). Therefore
\[
       \mathcal H_k
        =
        \mathcal M_k^+
        \oplus
        S_{2,k}
        \oplus\cdots\oplus
        S_{k,k}
        \oplus
        \ker B_k .
\]

Consequently every \(f\in\mathcal H_k\) has a unique decomposition as stated.
This proves the theorem.


\section{Applications}

The preceding results are structural: they identify the incremental and marginal spectral layers of the Green--Kubo covariance operator. Their concrete use, however, is model dependent, since the projections require conditional expectations under the stationary path law and, in the marginal layer, the action of the base transition operator \(K\). We illustrate this in two settings where these operations are explicit: reversible two-state Markov chains and Gaussian AR(1) chains.

\subsection{Reversible two-state chains}

We first consider the stationary reversible two-state Markov chain on
\(\{-1,1\}\) with transition matrix
\[
        K=
        \begin{pmatrix}
        p & q\\
        q & p
        \end{pmatrix},
        \qquad p+q=1.
\]
Its stationary distribution is uniform, and the nontrivial eigenvalue of \(K\) is
\[
        \rho=p-q.
\]
Thus, if \(x(s)=s\) denotes the coordinate function on \(\{-1,1\}\), then $Kx=\rho x$.

Equivalently,
\[
        E[X_{t+h}\mid X_t]=\rho^h X_t,
        \qquad
        E[X_tX_{t+h}]=\rho^h.
\]

It is useful to introduce the edge variables
\[
        R_i=X_iX_{i+1},\qquad i=1,\ldots,k-1.
\]
For the symmetric two-state chain, \(X_1,R_1,\ldots,R_{k-1}\) are independent,
with
\[
        P(R_i=1)=p,\qquad P(R_i=-1)=q,
        \qquad E[R_i]=\rho.
\]
Let
\[
        r_i=R_i-\rho.
\]
Then the variables \(r_i\) are centered and independent.  This gives an explicit
orthogonal coordinate system for functions of a \(k\)-block.

For \(k=2\), the incremental space is one-dimensional:
\[
        \mathcal I_2
        =
        \operatorname{span}\{R_1-\rho\}
        =
        \operatorname{span}\{X_1X_2-\rho\}.
\]
Indeed,
\[
        E[X_1X_2-\rho\mid X_1]=0,
        \qquad
        E[X_1X_2-\rho\mid X_2]=0.
\]
Thus the pure transition component at block length two is generated by the centered
edge variable \(r_1\).

The symmetry of the two-state chain also gives a natural parity decomposition. Under the global sign change \((x_1,\ldots,x_k)\mapsto (-x_1,\ldots,-x_k)\), a function is even if it is invariant and odd if it changes sign. In the edge coordinates, the variables \(R_i=X_iX_{i+1}\) are even, while \(X_1\) is odd. Thus even functions are functions of the edge variables alone, whereas odd functions are \(X_1\) times a function of the edge variables.

For \(k\ge3\), these edge coordinates give an explicit description of the even
and odd parts of \(\mathcal I_k\).  Since \(f\in\mathcal I_k\) must have zero
conditional expectation given both adjacent \((k-1)\)-blocks, the even part must
contain the two endpoint edge residuals \(r_1\) and \(r_{k-1}\).  Hence
\[
        \mathcal I_k^{\rm even}
        =
        r_1r_{k-1}\,
        \operatorname{span}
        \left\{
        \prod_{j\in J} r_j : J\subset\{2,\ldots,k-2\}
        \right\}.
\]
There is also an odd part, involving the centered two-sided innovation
\[
        X_1-\rho X_2.
\]
It is given by
\[
        \mathcal I_k^{\rm odd}
        =
        (X_1-\rho X_2)r_{k-1}\,
        \operatorname{span}
        \left\{
        \prod_{j\in J} r_j : J\subset\{2,\ldots,k-2\}
        \right\}.
\]
Thus
\[
        \mathcal I_k
        =
        \mathcal I_k^{\rm even}\oplus \mathcal I_k^{\rm odd},
        \qquad k\ge3.
\]

This representation gives more than the general theory alone. The general results identify \(\mathcal I_k\) abstractly as a two-sided innovation space, whereas the symmetry of the two-state chain supplies independent edge coordinates and hence explicit bases for its even and odd parts. The even component detects full-span dependence through the two endpoint edge residuals \(r_1\) and \(r_{k-1}\). The odd component is genuinely Markovian: it contains the boundary correction \(X_1-\rho X_2\), which is forced by the requirement of orthogonality to the right adjacent block.

The one-coordinate marginal spectrum is also explicit.  Since the only nonconstant eigenfunction
of the base transition operator is \(x\), with eigenvalue \(\rho\), Theorem~3.6 gives
the marginal eigenfunction
\[
        M_{\rho,k}(x_1,\ldots,x_k)
        =
        \frac{1}{1-\rho}x_1
        +
        \sum_{i=2}^{k-1}x_i
        +
        \frac{1}{1-\rho}x_k,
\]
with eigenvalue
\[
        \theta_k(\rho)=\frac{k-(k-2)\rho}{1-\rho}.
\]
Thus, for the reversible two-state chain, the full nonzero covariance structure
consists of the marginal eigenfunction above together with the shift-inherited
integer layers generated by the spaces \(S_{m,k}\), \(2\le m\le k\).

\subsection{Gaussian AR(1) chains}

We next consider a continuous-state example.  Let \(\{X_t\}\) be the stationary
Gaussian AR(1) chain
\[
        X_{t+1}=\rho X_t+\sqrt{1-\rho^2}\,\varepsilon_{t+1},
        \qquad |\rho|<1,
\]
where the variables \(\varepsilon_t\) are i.i.d. standard normal.  Then
\(X_t\sim N(0,1)\), the chain is reversible, and the transition operator \(K\)
acts on \(L^2_0(\gamma)\), where \(\gamma\) is the standard normal law.

The important feature of this example is that \(K\) has an explicit infinite
spectral decomposition.  If \(H_j\), \(j\ge1\), denote the normalized Hermite
polynomials, then
\[
        K H_j=\rho^j H_j.
\]
Thus the continuous state space produces an infinite marginal spectral layer:
each Hermite direction \(H_j\) gives a boundary-corrected marginal eigenfunction.
This contrasts with the two-state example, where the base transition operator has
only one nonconstant eigenfunction and hence only one marginal spectral direction.

The marginal layer of the \(k\)-block covariance operator is therefore completely
explicit.  For each \(j\ge1\), Theorem~\ref{T:marginal} gives the boundary-corrected marginal
eigenfunction
\[
        M_{j,k}(x_1,\ldots,x_k)
        =
        \frac{1}{1-\rho^j}H_j(x_1)
        +
        \sum_{i=2}^{k-1}H_j(x_i)
        +
        \frac{1}{1-\rho^j}H_j(x_k),
\]
with eigenvalue
\[
        \theta_k(\rho^j)
        =
        \frac{k-(k-2)\rho^j}{1-\rho^j}.
\]
Hence the one-coordinate marginal spectrum is the infinite family
\[
        \left\{\theta_k(\rho^j):j\ge1\right\}.
\]

The incremental part is also computable through Gaussian conditional
expectations.  Already for \(k=2\), the pure transition component of a two-variable
observable is obtained by subtracting its one-coordinate projections as in
Theorem~3.3.  For example, applying the \(k=2\) formula to
\[
        f(x,y)=xy-\rho
\]
gives
\[
        g(x,y)
        =
        xy-\rho
        -
        \frac{\rho}{1+\rho^2}
        \bigl[(x^2-1)+(y^2-1)\bigr],
\]
and this \(g\) belongs to \(\mathcal I_2\).  Thus even a simple quadratic transition
observable has a marginal part that must be removed before the
two-block incremental contribution is isolated.

For higher block lengths, the spaces \(S_{m,k}\) are obtained by embedding the
incremental components \(\mathcal I_m\) into all valid positions inside a \(k\)-block
and taking their block sums.  These give the shift-inherited integer eigenvalues
\(1,\ldots,k-1\), exactly as in the general theory.  The AR(1) example therefore
shows that the same decomposition applies beyond finite-state chains: the overlap
hierarchy remains discrete, while the marginal spectral layer reflects the infinite
Hermite spectrum of the continuous-state Markov operator.

\section{Conclusion}

We studied additive functionals of stationary Markov chains whose observables depend
on finite consecutive blocks.  The main result is that the intrinsic incremental space
\(\mathcal I_k\), defined as the part of a \(k\)-block observable orthogonal to the two
adjacent \((k-1)\)-block spaces, remains an eigenvalue-one subspace under Markov dependence.  For a stationary
Markov chain it admits the dynamic characterization  (\ref{E:intersecKernels}),
and hence it remains an eigenvalue-one subspace of the Green–Kubo covariance operator. Lower-block incremental components embed into longer blocks and produce the shift-inherited integer eigenvalues $1,\ldots, k-1$, while in the reversible case the remaining marginal layer is governed by the spectrum of the base transition operator. We also identify explicitly the transformed one-coordinate marginal spectrum and the corresponding boundary-corrected eigenfunctions.
The examples show how the decomposition becomes explicit in both finite-state and continuous-state settings, and suggest that the incremental space is of interest well beyond the present spectral computation.

First, the characterization in Theorem~\ref{T:universal} gives an operator-theoretic analogue of the classical problem of determining the effective block order, or finite-memory order, of a stationary process  (the memory order of a Markov chain of higher order): whether a process genuinely requires block length $k$ to be described, or whether every $k$-block observable is already representable through its two adjacent \((k-1)\)-blocks. 

Vanishing of \(\mathcal I_k\) is precisely the latter, in the closed \(L^2\) sense, and therefore indicates that the \(k\)-block observable structure reduces to adjacent \((k-1)\)-block information. In this way, the incremental space gives an operator-theoretic counterpart to likelihood-based and context-tree methods traditionally used for finite-memory processes~\cite{BuhlmannWyner1999,CsiszarTalata2006,Katz1981,Tong1975}. This also connects naturally to sliding-window pattern statistics
(Good~\cite{Good1953}, Pincus~\cite{Pincus1991}) and suggests possible applications to symbolic
sequence analysis, where window length selection is often a central practical issue.

Second, the two-sided nature of $\mathcal{I}_k$, orthogonal to information available from either the left or the right adjacent block--parallels the notion of predictive sufficiency developed in computational  mechanics, where a process's causal states (Shalizi and Crutchfield \cite{ShaliziCrutchfield2001}) capture exactly the information about the past relevant to predicting the future, no more and no less. The incremental space is a finite-block, Hilbert-space analogue of this idea: it isolates what a length-$k$ window contributes beyond what shorter windows already predict on both sides. Making this correspondence precise, and relating the spectral decomposition of \(B_k\)  to information-theoretic quantities such as the excess entropy, is a natural direction for future work.

More broadly, the present results suggest several further directions, including a systematic study of the rate at which \(\|f^{[k]}\|^2\)
decays with $k$ for structured processes, connections to predictive-state and state-compression methods for sequential data, and an asymptotic theory for order-selection procedures built directly from the spectral decomposition of \(B_k\).

\section*{Statements and Declarations}
The author declares no competing interests.



\begin{thebibliography}{99}

\bibitem{Alhakim2004} A. Alhakim, \textit{On the eigenvalues and
eigenvectors of an overlapping Markov chain}, Probability
Theory and Related Fields, 128 (2004), pp. 589-605.

\bibitem{AlhakimKawczakMolchanov2004} A. Alhakim, J. Kawczak, and S. Molchanov, \textit{On
the class of nilpotent Markov chains, I: the spectrum of
covariance operator}, Markov Processes and Related Fields
10 (2004), pp. 629-652.

\bibitem{AA2026} 
A.~Alhakim,
\newblock Overlapping window tests for correlation and trend,
\newblock Manuscript submitted for publication, 2026.


\bibitem{Brillinger1975}
D.~R. Brillinger,
\newblock Time Series: Data Analysis and Theory,
\newblock Holt, Rinehart and Winston, New York, 1975.

\bibitem{BuhlmannWyner1999}
P. B\"uhlmann and A. J. Wyner,
Variable length Markov chains,
\emph{The Annals of Statistics}, 27(2):480--513, 1999.

\bibitem{KLChung1982}
K. L. Chung,
\textit{Lectures from Markov Processes to Brownian Motion},
Grundlehren der mathematischen Wissenschaften, Vol. 249,
Springer-Verlag, New York, 1982.

\bibitem{CsiszarTalata2006}
I. Csisz\'ar and Zs. Talata,
Context tree estimation for not necessarily finite memory processes, via BIC and MDL,
\emph{IEEE Transactions on Information Theory}, 52(3):1007--1016, 2006.

\bibitem{DehlingWendler2010}
H.~Dehling and M.~Wendler,
\newblock Central limit theorem and the bootstrap for \(U\)-statistics of strongly
mixing data,
\newblock Journal of Multivariate Analysis \textbf{101} (2010), no.~1, 126--137.

\bibitem{Good1953} I. J. Good, ``\textit{The serial test for sampling numbers and other tests of randomness}'',
In Proceedings of Cambridge Philosophical Society,  49, (1953), 276-284.

\bibitem{Gordin1969}
M.~I. Gordin,
\newblock The central limit theorem for stationary processes,
\newblock Soviet Math. Dokl. \textbf{10} (1969), 1174--1176.
\newblock [Translated from Dokl. Akad. Nauk SSSR \textbf{188} (1969), 739--741.]

\bibitem{HoHsing1997}
H.-C. Ho and T.~Hsing,
\newblock Limit theorems for functionals of moving averages,
\newblock Ann. Probab. \textbf{25} (1997), no.~4, 1636--1669.

\bibitem{Hoeffding1948}
W.~Hoeffding,
\newblock A class of statistics with asymptotically normal distribution,
\newblock Ann. Math. Statist. \textbf{19} (1948), 293--325.

\bibitem{IlIdrissi2025}
M.~Il Idrissi, N.~Bousquet, F.~Gamboa, B.~Iooss, and J.-M.~Loubes,
\newblock Hoeffding decomposition of functions of random dependent variables,
\newblock Journal of Multivariate Analysis \textbf{208} (2025), 105444.

\bibitem{KacMurdockSzego1953}
M.~Kac, W.~L. Murdock, and G.~Szegő,
\newblock On the eigen-values of certain Hermitian forms,
\newblock Journal of Rational Mechanics and Analysis \textbf{2} (1953), 767--800.

\bibitem{Katz1981}
R. W. Katz,
On some criteria for estimating the order of a Markov chain,
\emph{Technometrics}, 23(3):243--249, 1981.

\bibitem{KipnisVaradhan1986}
C.~Kipnis and S.~R.~S. Varadhan,
\newblock Central limit theorem for additive functionals of reversible Markov processes
and applications to simple exclusions,
\newblock Comm. Math. Phys. \textbf{104} (1986), no.~1, 1--19.


\bibitem{Pincus1991} S. M. Pincus, \textit{Approximate entropy as a measure of system complexity},  Proceedings of the National Academy of Sciences, USA,  88, (1991), pp. 2297–2301.


\bibitem{ShaliziCrutchfield2001}
C. R. Shalizi and J. P. Crutchfield,
Computational mechanics: pattern and prediction, structure and simplicity,
\emph{Journal of Statistical Physics}, 104:817--879, 2001.


\bibitem{Tong1975}
H. Tong,
Determination of the order of a Markov chain by Akaike's information criterion,
\emph{Journal of Applied Probability}, 12:488--497, 1975.




\bibitem{Wu2005}
W.~B. Wu,
\newblock Nonlinear system theory: another look at dependence,
\newblock Proc. Natl. Acad. Sci. USA \textbf{102} (2005), no.~40, 14150--14154.



\end{thebibliography}
\end{document}